\documentclass[11pt]{article}
\usepackage{latexsym, amssymb, amsmath}

\def \F {{\mathbb F}}

\def \Tr {{\rm Tr_n}}
\def \T {{\rm Tr}}
\def \la {\langle}
\def \ra {\rangle}
\newtheorem{theorem}{Theorem}

\newtheorem{lemma}{Lemma}
\newtheorem{proposition}{Proposition}
\newtheorem{remark}{Remark}
\newtheorem{example}{Example}
\newtheorem{corollary}{Corollary}
\def\bproof{{\it Proof}: }
\def\eproof{\hfill$\Box$}

\begin{document}

\begin{center}
{\Large \bf (Not) weakly regular univariate bent functions}
\end{center}
\begin{center}
Ay\c ca \c Ce\c smelio\u glu$^1$, Wilfried Meidl$^2$, \\[.7em]
$^1$ Department of Mathematics, Otto-von-Guericke-University, 39106 Magdeburg, Germany. \\
email: cesmelioglu@gmail.com \\
$^2$ Sabanc\i\ University, MDBF, Orhanl\i, Tuzla, 34956 \.Istanbul, Turkey. \\
email: wmeidl@sabanciuniv.edu
\end{center}

\begin{abstract}
In this article a procedure to construct bent functions from $\F_{p^n}$ to
$\F_p$ by merging plateaued functions which are bent on ($n-2$)-dimensional 
subspaces of $\F_{p^n}$ is presented. Taking advantage of such classes
of plateaued functions with a simple representation as monomials and
binomials, we obtain infinite classes of bent functions with a fairly
simple representation. In particular we present the first direct 
construction of univariate not weakly regular bent functions, and give 
one class explicitly in a simple representation with binomials. 
%
\end{abstract}

{\it Keywords:}
Bent function; partially bent function; Fourier transform;
not weakly regular; quadratic function; polynomial.

\section{Introduction}

Let $V_n$ be an $n$-dimensional vector space over the prime field $\F_p$.
A function $f:V_n\rightarrow\F_p$ is called a {\it bent function} if its
{\it Fourier transform} $\widehat{f}$ defined by
\[ \widehat{f}(b) = \sum_{x\in V_n}\epsilon_p^{f(x)-<b,x>} \]
satisfies $|\widehat{f}(b)|^2 = p^n$ for all $b \in V_n$, where
$\epsilon_p = e^{2\pi i/p}$ and $<,>$ denotes any (non-degenerate)
inner product on $V_n$.
Classical representations for bent functions are the multivariate representation
where $V_n = \F_p^n$, in this case one may use the conventional dot product as inner
product, and the univariate representation where $V_n = \F_{p^n}$, in which case
one may use $<b,x> = \Tr(bx)$ as inner product, where $\Tr(z)$
denotes the absolute trace of $z\in\F_{p^n}$. 

For $p=2$, bent functions can only exist when $n$ is even, the {\it Fourier coefficients}
$\widehat{f}(b)$ are then obviously $\pm 2^{n/2}$. For $p > 2$ bent functions exist for both, $n$
even and $n$ odd. For the Fourier coefficients we then always have (cf. \cite{hk})
\begin{equation}
\label{(2)} p^{-n/2}\widehat{f}(b) =
\left\{\begin{array}{r@{\quad:\quad}l}
\pm \epsilon_p^{f^*(b)} & n\;\mbox{even or}\;n\;\mbox{odd and}\;p \equiv 1\bmod 4 \\
\pm i\epsilon_p^{f^*(b)} & n\;\mbox{odd and}\;p \equiv 3\bmod 4
\end{array}\right.
\end{equation}
where $f^*$ is a function from $V_n$ to $\F_p$.
A bent function $f:V_n\rightarrow\F_p$ is called {\it regular} if for all $b \in V_n$
\[ p^{-n/2}\widehat{f}(b) =\epsilon_p^{f^*(b)}. \]
When $p=2$, a bent function is trivially regular, and as can be seen from $(\ref{(2)})$, for
$p>2$ a regular bent function can only exist for even $n$, and for odd $n$ when $p\equiv 1\bmod 4$.
A function $f:V_n\rightarrow\F_p$ is called {\it weakly regular} if, for all $b \in \F_p^n$, we have
\[ p^{-n/2}\widehat{f}(b) =\zeta\ \epsilon_p^{f^*(b)} \]
for some complex number $\zeta$ with $|\zeta|=1$, otherwise it is called {\it not weakly regular}.
By (\ref{(2)}), $\zeta$ can only be $\pm 1$ or $\pm i$. Note that regular implies weakly regular.
All classical construction of bent functions yield (weakly) regular bent functions.

Based on earlier constructions of Boolean bent functions from near-bent functions
(see \cite{cpt,lg}) in \cite{agw,aw1} constructions of $p$-ary bent functions have been 
presented. In these constructions functions in lower dimensions are merged to a bent function 
by adjoining variables. The resulting bent functions are given in multivariate form \cite{agw,aw}, 
or as functions from $\F_{p^{n-1}}\times\F_p$ to $\F_p$ \cite{aw1}. The construction turns out to
be very powerful, for instance the first infinite classes of not weakly regular bent functions 
have been obtained. Until then only sporadic examples of not weakly regular bent 
functions were known, all given as univariate polynomials of the form $f(x) =
\Tr(g(x))$ for a polynomial $g(x)\in\F_{p^n}[x]$, see \cite{hk,hk1,tyz}. 

The objective of this paper is to develop an equivalent construction for the
univariate case, i.e. for functions from $\F_{p^n}$ to $\F_p$, which is more involved
as we cannot simply add variables to the finite field $\F_{p^n}$. Amongst others, sets
of functions which are bent on ($n-2$)-dimensional subspaces of $\F_{p^n}$ are required.
We take advantage of particularly simple representations of some classes of such functions
(as monomials and binomials), and thereby obtain bent functions from $\F_{p^n}$ to $\F_p$
in a simple representation. Among those, we present the first direct construction of infinite 
classes of univariate not weakly regular bent functions. 

In Section \ref{sec2} we develop the principles of the construction and we give an explicit 
formula for bent functions from $\F_{p^n}$ to $\F_p$ obtained from functions which are bent
on subspaces. This description of the functions enables us also to find the representation 
with a unique polynomial of degree at most $p^n-1$. In Section \ref{sec3} we use classes
of monomials and binomials to obtain univariate (not) weakly regular bent functions in a
simple representation. Explicitly we describe an infinite class of not weakly regular bent
functions. Some examples are given in the appendix.

\section{A construction of bent polynomials}
\label{sec2}

For a function $f:V_n\rightarrow\F_p$ and an element $a\in V_n$, the {\it derivative} $D_af:V_n\rightarrow\F_p$
of $f$ in direction $a$ is defined by $D_af(x) = f(x+a)-f(x)$. As well known, $f$ is bent if and only 
if $D_af$ is balanced for all nonzero $a\in V_n$, see \cite{n}. An element $a\in V_n$ for which
$D_af$ is constant is called a {\it linear structure} of $f$. As easily seen, the set $\Lambda$ of the
linear structures of $f$ forms a subspace of $V_n$, which we call the {\it linear space} of the function $f$.
We have
\begin{equation}
\label{linear}
f(x+a) = f(x)+f(a)-f(0)\;\mbox{for all}\;a\in \Lambda, x\in V_n.
\end{equation}
In particular, if $f(0) = 0$, then equation $(\ref{linear})$ implies that $f$ is linear on $\Lambda$.

A function $f:V_n\rightarrow\F_p$ is called {\it partially bent} if
for all $a\in V_n$ the derivative $D_af$ is either balanced or
constant. The set of partially bent functions is a subset of the set
of {\it plateaued functions}, which is the set of functions
$f:V_n\rightarrow\F_p$ for which $\widehat{f}(b) = 0$ or
$|\widehat{f}(b)| = p^{(n+s)/2}$ for all $b\in V_n$ and a fixed
integer $s$, $0\le s\le n$, depending on $f$. This can easily be
seen in the calculations below, applying the standard Welch-squaring
method. In accordance with \cite{agw}, we call plateaued functions
$f$ from $V_n$ to $\F_p$ for which $s=1$ as {\it near-bent functions}. We
remark that when $p=2$, also the term {\it semi-bent function} is used for
plateaued functions
with $s=1$ and $n$ odd or $s=2$ and $n$ even, see \cite{cpt}. \\
Let $f$ be a partially bent function with linear space $\Lambda$ of dimension $s$ as a
subspace of $V_n$. Without loss of generality we will always suppose that $f(0)=0$, and hence
$f(x+a) = f(x) + f(a)$ if $a$ is a linear structure of $f$. We then have
\begin{eqnarray*}
|\widehat{f}(b)|^2 & = &
\sum_{x,y\in V_n}\epsilon_p^{f(x)-f(y)-<b,x-y>} =
\sum_{y,z\in V_n}\epsilon_p^{f(y+z)-f(y)-<b,z>} \\
& = & \sum_{z\in V_n}\epsilon_p^{f(z)-<b,z>}
\sum_{y\in V_n}\epsilon_p^{f(y+z)-f(y)-f(z)}.
\end{eqnarray*}
Using that $f(y+z)-f(y)-f(z)$ is balanced as a function in variable $y$ if $z\not\in \Lambda$, we get
\begin{equation}
\label{support}
|\widehat{f}(b)|^2 = p^n\sum_{z\in \Lambda}\epsilon_p^{f(z)-<b,z>}
= \left\{\begin{array}{lr}
p^{n+s} & \mbox{if}\;f(z) - <b,z> \equiv 0\;\mbox{on}\; \Lambda, \\
0 & \mbox{otherwise,}\end{array}\right.
\end{equation}
where in the last step we use that $f$ is linear on $\Lambda$.
Defining the support of the Fourier transform of $f$ by
$supp(\widehat{f}) := \{b\in V_n\;:\;\widehat{f}(b) \ne 0\}$, it follows from
equation $(\ref{support})$ that $b\in supp(\widehat{f})$ if and only if $f(z)-<b,z> \equiv 0$ on $\Lambda$.
{\it Parseval's identity}
\[ \sum_{b\in V_n}|\widehat{f}(b)|^2 = p^{2n}, \]
then implies $|supp(\widehat{f})| = p^{n-s}$.
\begin{remark}
As it can be seen from equation $(\ref{support})$, $supp(\widehat{f})$ depends on
the inner product $<,>$ which is used. Consequently, to be precise one may define the support
of $\widehat{f}$ with respect to the inner product $<,>$. As well known the absolute values
appearing in the Fourier spectrum $\{\widehat{f}(b)\;|\;b\in V_n\}$ of $f$ are independent of
the (non-degenerate) inner product. In particular, the property of being $s$-plateaued is independent
from $<,>$.
\end{remark}
We will use the following result on partially bent functions (see \cite[Theorem]{c}).
\begin{lemma}
\label{Carlemma}
Let $f:V_n\rightarrow\F_p$ be a partially bent function with linear space $\Lambda$ and let
$\Lambda^c$ be any complement of $\Lambda$ in $V_n$. Then $f$ restricted to $\Lambda^c$ is
a bent function.
\end{lemma}
A well understood class of partially bent functions is the class of quadratic functions,
see \cite{agw,hk}. For more information on partially bent functions we refer the reader to \cite{c}.

For the construction of bent functions from $\F_{p^n}$ to $\F_p$ we will employ partially bent
functions $f$ with a two dimensional linear space $\Lambda = \la\beta_1,\beta_2\ra$.
For simplicity we fix an inner product $<u,v> = \Tr(\delta uv)$ on $\F_{p^n}$ with respect to which the orthogonal
complement $\Lambda^\perp$ of $\Lambda$ is a complement $\Lambda^c$ of $\Lambda$, and we further
suppose that $<\beta_1,\beta_2> = 0$. Note that this implies $<\beta_1,\beta_1> = \ell \ne 0$ and
$<\beta_2,\beta_2> = t \ne 0$. Conversely, the properties $<\beta_1,\beta_2> = 0$, $<\beta_1,\beta_1> = \ell \ne 0$ and
$<\beta_2,\beta_2> = t \ne 0$ imply that $\Lambda\cap\Lambda^\perp = \{0\}$, and hence
that $\Lambda^\perp$ is a complement of $\Lambda$.
We remark that given $\beta_1,\beta_2 \in \F_{p^n}$ (linearly independent over $\F_p$), one can always find $\delta \in\F_{p^n}$ 
such that $\Tr(\delta\beta_1\beta_2) = 0$, $\Tr(\delta\beta_1^2) \ne 0$ and $\Tr(\delta\beta_2^2) \ne 0$.
Some properties of this inner product are used in the proof of Lemma \ref{near-bent}, 
some are used in the proof of Theorem \ref{theo} below, see also Remark \ref{inepro}.

For a partially bent function $f:\F_{p^n}\rightarrow\F_p$
with linear space $\Lambda = \la\beta_1,\beta_2\ra$ we fix the following notation:
\begin{itemize}
\item[-] $f_{|n-2}$ is the function $f$ restricted to $V_{n-2}:=\Lambda^c$,
\item[-] $f_{|n-1}$ is the function $f$ restricted to $V_{n-1}:=\la\beta_2\ra + \Lambda^c$,
\item[-] $\widehat{f}_{|n-2}$ is the Fourier transform of $f_{|n-2}$,
\item[-] $\widehat{f}_{|n-1}$ is the Fourier transform of $f_{|n-1}$.
\end{itemize}
The following lemma shows that $f$ given as above 
is near-bent with linear
space $\la\beta_2\ra$ when restricted to $V_{n-1}=\la\beta_2\ra + \Lambda^c$. We consider the case
of $\Lambda^c = \Lambda^\perp$ with an inner product defined as above.
\begin{lemma}
\label{near-bent}
Let $f:\F_{p^n}\rightarrow\F_p$ be a partially bent function with linear space
$\Lambda = \la\beta_1,\beta_2\ra$, and let
$V_{n-1} = \la\beta_2\ra + \Lambda^c$. Then $f|_{n-1}:V_{n-1}\rightarrow\F_p$ is near-bent
with linear structure $\la\beta_2\ra$.
\end{lemma}
\bproof
Supposing that $<,>$ satisfies the above described conditions, we observe that $<,>$ also defines
a non-degenerate inner product on $V_{n-1}$ and on $V_{n-2} = \Lambda^\perp = \Lambda^c$. Further we know that
$<\beta_2,\beta_2> = t \ne 0$. Obviously the linear structure $\beta_2$ of $f$ is also a linear
structure of $f_{|n-1}$. For an element $b = b_2\beta_2+y\in V_{n-1}$, $b_2\in\F_p$,
$y\in V_{n-2} = \Lambda^c$, with $(\ref{linear})$ we then have
\begin{eqnarray*}
\widehat{f}_{|n-1}(b) & = &
\sum_{c\in\F_p}\sum_{x\in\Lambda^c}\epsilon_p^{f(c\beta_2+x)-<b_2\beta_2+y,c\beta_2+x>}
= \sum_{c\in\F_p}\epsilon_p^{f(c\beta_2)-b_2ct}\sum_{x\in\Lambda^c}\epsilon_p^{f(x)-<x,y>} \\
& = & \widehat{f}_{|n-2}(y)\sum_{c\in\F_p}\epsilon_p^{cf(\beta_2)-b_2ct} =
\left\{\begin{array}{ll}
0 & \mbox{if}\;b_2 \ne f(\beta_2)/t, \\
p\widehat{f}_{|n-2}(y) & \mbox{if}\;b_2 = f(\beta_2)/t.
\end{array}
\right.
\end{eqnarray*}
Consequently by Lemma \ref{Carlemma}, $f_{|n-1}$ is near-bent.
\eproof\\[.5em]
With the following proposition we can generate sets of $p$ near-bent functions on $V_{n-1}$ such that every
element of $V_{n-1}$ is in the support of the Fourier transform for exactly one function in the set.
\begin{proposition}
\label{prop}
Let $f_0,\ldots,f_{p-1}$ be partially bent functions from $\F_{p^n}$ to $\F_p$, all with the same
linear space $\Lambda = \langle \beta_1,\beta_2 \rangle$, and
for $0\le k\le p-1$ let $\gamma_k\in V_{n-1}$ be such that
\begin{equation}
\label{gammas}
f_k(\beta_2) +<\gamma_k,\beta_2> = f_0(\beta_2)+k.
\end{equation}
The functions $g_{k_{|n-1}}:V_{n-1}\rightarrow\F_p$ defined by
\begin{equation}
\label{gk}
g_k(x) = f_k(x)+<\gamma_k,x>,\quad k=0,1,\ldots,p-1
\end{equation}
(restricted to $V_{n-1}$) form a set of near-bent functions with $supp(\widehat{g_j}_{|n-1})\cap supp(\widehat{g_k}_{|n-1})
= \emptyset$ if $j\ne k$.
\end{proposition}
\bproof
By the above discussion it is guaranteed that $f_{k_{|n-1}}$, i.e. $f_k$ restricted to $V_{n-1}$,
is near-bent with linear space $\la\beta_2\ra$ for all $k=0,\ldots,p-1$. We have to show that
the addition of the linear functions $<\gamma_k,x>$ to the functions $f_k$ separates the supports
of the Fourier transforms of the corresponding functions on $V_{n-1}$. By
equation $(\ref{support})$, $b\in V_{n-1}$ is an element of $supp(\widehat{g_k}_{|n-1})$
(with respect to the inner product $<,>$) if and only if $g_k(\beta_2)-<b,\beta_2> = 0$.
Suppose that $b \in supp(\widehat{g_k}_{|n-1}) \cap supp(\widehat{g_j}_{|n-1})$
for some $0\le k,j\le p-1$. Then
\begin{eqnarray*}
0 & = & f_k(\beta_2)+<\gamma_k,\beta_2>-<b,\beta_2> = f_0(\beta_2)+k-<b,\beta_2> \\
& = & f_j(\beta_2)+<\gamma_j,\beta_2>-<b,\beta_2> = f_0(\beta_2)+j-<b,\beta_2>,
\end{eqnarray*}
and hence $j=k$. Since as a consequence of Parseval's identity we have
$|supp(\widehat{g_k}_{|n-1})| = p^{n-2}$, the result follows.
\eproof\\[.5em]
We remark that though the support of a near-bent function depends on the considered inner product,
by the method of Proposition \ref{prop} one obtains a set of $p$ near-bent functions such that every
element $b$ of $V_{n-1}$ is in the support of the Fourier transform for exactly one function,
independent of the inner product used.
\begin{theorem}
\label{theo}
Let $<,>$ be an inner product on $\F_{p^n}$ and let $\beta_1\in\F_{p^n}$ be
such that $<\beta_1,\beta_1> = \ell \ne 0$, let $g_0,\ldots,g_{p-1}$ be functions from $\F_{p^n}$ to
$\F_p$ with linear structure $\beta_1$, and suppose that the restrictions $g_{k_{|n-1}}$ to
$V_{n-1} = \la\beta_1\ra^\perp$ are near-bent and $supp(\widehat{g_j}_{|n-1})\cap
supp(\widehat{g_k}_{|n-1}) =\emptyset$ if $j\ne k$. With $\gamma = \ell^{-1}\beta_1$, the
function $F:\F_{p^n}\rightarrow\F_p$
\begin{equation}
\label{F}
F(x) = -\sum_{k=0}^{p-1}\prod_{j=0\atop j\ne k}^{p-1}(<\gamma,x>-j)g_k(x)
\end{equation}
is bent.
\end{theorem}
\bproof
Since $<\beta_1,\beta_1> = \ell \ne 0$, the orthogonal complement $V_{n-1}=\la\beta_1\ra^\perp$ is a
complement of $\la\beta_1\ra$, and $<,>$ is a non-degenerate inner product on $V_{n-1}$.
Let $x=c_1\beta_1+y$, $c_1\in\F_p$, $y\in V_{n-1}$, be the unique representation of $x\in\F_{p^n}$
as a sum of elements of $V_{n-1}^\perp = \la\beta_1\ra$ and $V_{n-1}$. Then using
$<\gamma,\beta_1> = \ell^{-1}<\beta_1,\beta_1> = 1$ and $\la\beta_1\ra = V_{n-1}^\perp$ we obtain
\[ <\gamma,x>-j = c_1<\gamma,\beta_1> + <\gamma,y> -j = c_1-j. \]
Hence $\prod_{j=0\atop j\ne k}^{p-1}(<\gamma,x>-j)$ vanishes if $k \ne c_1$ and for $k = c_1$
we have $\prod_{j=0\atop j\ne k}^{p-1}(<\gamma,x>-j) = -1$. Consequently,
\[ F(x) = g_{c_1}(x) = g_{c_1}(c_1\beta_1 + y) = g_{c_1}(c_1\beta_1) + g_{c_1}(y), \]
where in the last step we use that $\beta_1$ is a linear structure of $g_{c_1}$. \\
Let $b = b_1\beta_1 + z$, $b_1\in\F_p$, $z\in V_{n-1}$. Again using $V_{n-1} = \langle\beta_1\rangle^\perp$
we then get
\begin{eqnarray*}
\widehat{F}(b_1\beta_1+z) & = & \sum_{y\in V_{n-1}\atop c\in\F_p}\epsilon_p^{g_c(c\beta_1+y)-
<b_1\beta_1+z,c\beta_1+y>} \\
& = & \sum_{c\in\F_p}\epsilon_p^{g_c(c\beta_1)-b_1c\ell}\sum_{y\in V_{n-1}}\epsilon_p^{g_c(y)-\Tr(zy)}
= \sum_{c\in\F_p}\epsilon_p^{g_c(c\beta_1)-b_1cl}\widehat{g_c}_{|n-1}(z).
\end{eqnarray*}
Since every $z\in V_{n-1}$ is in the support of $\widehat{g_c}_{|n-1}$ of exactly one $c$,
for this $c$ we then have
\[ \widehat{F}(b_1\beta_1+z) = \epsilon_p^{g_c(c\beta_1)-b_1c\ell}\widehat{g_c}_{|n-1}(z), \]
and therefore $|\widehat{F}(b_1\beta_1+z)| = p^{n/2}$
\eproof\\[.5em]
The near-bent functions on $V_{n-1}$ constructed in Proposition \ref{prop} together with the inner
product $<,>$ considered in the proposition satisfy the assumptions of Theorem \ref{theo}, and we
can suggest the following procedure for constructing bent polynomials:
\begin{itemize}
\item[-] Choose $p$ partially bent functions $f_k:\F_{p^n}\rightarrow\F_p$, $0\le k\le p-1$,
all with the same $2$-dimensional linear space $\Lambda = \langle\beta_1,\beta_2\rangle$.
\item[-] Choose an inner product $<u,v> = \Tr(\delta uv)$ on $\F_{p^n}$ such that \linebreak
$<\beta_1,\beta_2> = 0$, and $\Lambda^\perp$ with respect
to this inner product is a complement of $\Lambda$. We remark that therefore it is sufficient that
$\delta$ satisfies $\Tr(\delta\beta_1\beta_2) = 0$, $\Tr(\delta\beta_1^2) = \ell \ne 0$ and
$\Tr(\delta\beta_2^2) = t \ne 0$.
\item[-] For $k = 0,\ldots,p-1$ choose $\gamma_k\in V_{n-1} = \la\beta_1\ra^\perp$ which satisfy
equation $(\ref{gammas})$, to obtain the functions $g_k$ defined as in equation $(\ref{gk})$.
We emphasize that such elements $\gamma_k$ always exist.
\item[-] With $\gamma = \ell^{-1}\beta_1$ construct the function $F:\F_{p^n}\rightarrow\F_p$ given
as in equation $(\ref{F})$.
\end{itemize}
\begin{remark}
\label{inepro} For the construction of bent functions it is
sufficient to assume that $<\beta_1,\beta_1>=\ell\ne 0$ and
$<\beta_1,\beta_2>=0$. For $V_{n-1}=\la\beta_2\ra + \Lambda^c$ in
Theorem \ref{theo} one can then take $\la\beta_1\ra^\perp$. In this case, all
properties needed in the proof of Theorem \ref{theo} hold but $\Lambda^\perp$ may not be a complement of
$\Lambda$ with respect to this inner product. Then a different inner
product has to be considered for the proof of Lemma \ref{near-bent}
where the orthogonality of $\beta_2$ and $\Lambda^c$, and
$<\beta_2,\beta_2> = t \ne 0$ is needed.
\end{remark}
The representation of the bent function in Theorem \ref{theo} with a closed formula enables also 
the determination of the corresponding unique polynomial in $\F_{p^n}[x]$ of degree at most $p^n-1$.
As one may expect, in general this representation does not look simple at all (see the examples in 
Section \ref{sec3} and in the appendix). As pointed out in the following corollary, we have a rather 
simple representation for the functions defined in ($\ref{F}$) originated in their construction principle.
\begin{corollary}
\label{simple}
Let $f_k:\F_{p^n}\rightarrow\F_p$, $0\le k\le p-1$, be partially bent functions all with the same 
$2$-dimensional linear space $\Lambda = \langle\beta_1,\beta_2\rangle$ and let $< u,v > = \Tr(\delta uv)$ 
for an element $\delta \in \F_{p^n}$ such that $< \beta_1,\beta_1 > = l \ne 0$, and \linebreak
$< \beta_1,\beta_2 > = 0$. For $k = 0,\ldots,p-1$, let 
$\gamma_k\in V_{n-1} = \la\beta_1\ra^\perp$ be such that \linebreak
$f_k(\beta_2) + < \gamma_k,\beta_2 > = f_0(\beta_2) + k$.
Then with $g_k(x) := f_k(x) + < \gamma_k,x >$, $k = 0,\ldots,p-1$, and $\gamma = l^{-1}\beta_1$ the function
\[ F(x) = g_{<\gamma,x>}(x) \]
is a bent function.
\end{corollary}
\bproof
Observing that all requirements for the function defined as in $(\ref{F})$ to be bent are satisfied, the statement
follows since $(\ref{F})$ reduces to $g_k(x)$ if $< \gamma, x > = k$.
\eproof

\section{(Not) weakly regular bent polynomials from quadratic monomials and binomials}
\label{sec3}

The simplicity of the representation of the bent functions in Corollary \ref{simple}
equals the simplicity of the ingredient partially bent functions $f_k$. Hence the
above procedure for the construction of bent functions motivates the study of partially 
bent functions from $\F_{p^n}$ to $\F_p$ with a $2$-dimensional linear space and a simple 
representation. Since all quadratic functions are partially bent, it is natural to analyse 
elements of this class of functions (see e.g. \cite{cpt,mt}) starting with monomials and binomials (in trace form).

For quadratic monomials $f(x) = \Tr(\alpha x^{p^r+1})$ the Fourier 
coefficients are known. The subsequent lemma gives the conditions under which we
have a linear space of dimension $2$. For a proof for $p$ odd we
refer to \cite[Theorem 1]{aw1}, the case $p=2$ follows
straightforward with the same approach.
\begin{lemma}
\label{monomial}
Let $f(x) = \Tr(\alpha x^{p^r+1})$, let $g$ be a primitive element
of $\F_{p^n}$, and suppose that $\alpha = g^c$.
\begin{enumerate}
\item If $p=2$, then $f$ is $2$-plateaued if and only if
\begin{itemize}
\item[(i)] $n\equiv 0\bmod 4$, $r$ is odd and $3$ divides $c$, or
\item[(ii)] $n\equiv 2\bmod 4$, and $r$ is even or $3$ divides $c$.
\end{itemize}
\item If $p$ is odd, then $f$ is $2$-plateaued if and only if $n$ is even, $r$ is odd, and
$c$ satisfies the equation $y(p^{2}-1)+c(p^r-1) = (p^n-1)/2$ for some integer $y$.
\end{enumerate}
\end{lemma}
The following proposition presents an infinite class of quadratic binomials with a $2$-dimensional 
linear space.
\begin{proposition}\label{binomial}
Let $p$ be an odd prime, $n$ a positive integer divisible by $3$, and let $r=2^\kappa$ for an integer $\kappa\ge 0$.
The quadratic binomial $f(x) = \Tr(\frac{p+1}{2}x^2+x^{p^r+1})$ has a $2$-dimensional linear space
if and only if $\kappa = 0$, or $\kappa \ge 1$ and $n$ is odd.
\end{proposition}
\bproof
With the standard Welch-squaring method we see that the linear space $\Lambda$ of $f$ is the
kernel of the linearized polynomial (cf. \cite[Equation(3.2)]{agw})
\[ L(x)=x^{p^r}+x^{p^{2r}}+x \in \F_{p^n}[x]. \]
The dimension of $\Lambda$ is then the degree of $\gcd(x^n-1,A(x))=a(x)$ where
$A(x)=x^{r}+x^{2r}+1$ is the associate of $L(x)$, if $a(x) = \sum a_ix^i$, then the kernel of
$L(x)$ is the set of all solutions of $\sum a_ix^{p^i}$, see \cite[p.118]{ln}.
We recall that when $n=n_1p^v$, $\gcd(n_1,p) = 1$, then the polynomial $x^n-1 \in \F_{p^n}[x]$
can be factored as
\[ x^n-1=\prod_{d|n_1}(\Phi_d(x))^{p^v}, \]
where $\Phi_d(x)\in \F_{p^n}[x]$ is the $d$th cyclotomic polynomial which has degree $\varphi(d)$.
Using that $\Phi_6(x)=x^2-x+1$ and
$\Phi_6(x^{2^j})=\Phi_{6\cdot 2^j}(x)$ for $j \geq 0$, we see that for $r = 2^\kappa$, $\kappa\ge 0$,
the polynomial $A(x)$ factors as
\begin{align*}
A(x) &= x^{2^{\kappa+1}}+x^{2^\kappa}+1 = (x^2+x+1)\prod_{j=1}^{\kappa}(x^{2^j}-x^{2^{j-1}}+1)\\
&=\Phi_3(x)\prod_{j=0}^{\kappa-1}\Phi_6(x^{2^j})=\Phi_3(x)\prod_{j=0}^{\kappa-1}\Phi_{6\cdot 2^j}(x)
= \prod_{j=0}^\kappa\Phi_{3\cdot 2^j}(x).
\end{align*}
(For $\kappa=0$ the empty product $\prod_{j=1}^{\kappa}(x^{2^j}-x^{2^{j-1}}+1)$ is defined as $1$.)
Trivially the statement of the proposition holds for $\kappa=0$. If $p \ne 3$ then for all
cyclotomic polynomials $\Phi_m(x)$ which appear in the above factorizations of $x^n-1$ and of
$A(x)$ we have $\gcd(m,p) = 1$. Using that under this condition different cyclotomic polynomials are
relatively prime, for $\kappa\ge 1$ we have $\gcd(x^n-1,A(x)) = \Phi_3(x)$ when
$n$ is odd, and $\Phi_3(x)\Phi_6(x)$ divides $\gcd(x^n-1,A(x))$ when $n$ is even.
Since $\Phi_3(x)=(x+2)^2$ for $p=3$, we have $x^n-1 = (x-1)\Phi_3(x)\prod_{d|n_1\atop d>1}(\Phi_d(x))^{3^v}$
and $A(x) = \Phi_3(x)\prod_{j=0}^{\kappa-1}(x^{2^j}+1)^2$. The roots of $x^{2^j}+1$ are primitive
$2^{j+1}$th roots of unity, which are not roots of $x^n-1$ when $n$ is odd. If $n$ is even,
then obviously $\Phi_3(x)\Phi_2(x) = (x-1)^2(x+1)$ divides $\gcd(x^n-1,A(x))$.
\eproof

With Lemma \ref{monomial} and Proposition \ref{binomial} infinitely many bent functions given as in Theorem \ref{theo} and 
Corollary \ref{simple} can be constructed, which have a fairly simple description using monomials (plus a linear term) and binomials, respectively. 
We first give an example for the construction of a Boolean bent function with a simple representation using monomials.
Some further examples are given in the appendix.
\begin{example}
Let $g$ be a primitive element of $\F_{2^6}$ satisfying
$g^6=g^4+g^3+g+1$, i.e. $g$ is a root of the primitive polynomial
$x^6+x^4+x^3+x+1\in\F_2[x]$. According to Lemma \ref{monomial}, the
functions
\[ f_0(x)=\T_6(gx^5)\;\mbox{and}\;f_1(x)=\T_6(g^{22}x^5) \]
are $2$-plateaued functions from $\F_{p^n}$ to $\F_p$. For both functions the linear space is
$\Lambda=\la\beta_1,\beta_2\ra$ with $\beta_1=g^{25}=g^3+g^2+1,\beta_2=g^{46}=g^4+g^3+g^2+1$.
We then have $\T_6(\beta_1^2)=1$, $\T_6(\beta_1\beta_2)=0$, $\T_6(\beta_2^2)\ne 0$.
We apply Proposition \ref{prop} choosing $\gamma_0=0$ and $\gamma_1=g^3$, and obtain
the near-bent functions from $V_5$ to $\F_2$
\[ g_{0_{|5}}(x)=\T_6(gx^5)\;\mbox{and}\;g_{1_{|5}}(x)=\T_6(g^2x^5+g^3x) \]
with $supp(\widehat{g_0}_{|5})\cap supp(\widehat{g_1}_{|5})=\emptyset$.
We apply Theorem \ref{theo} with $\gamma=\beta_1=g^{25}$ and obtain the
bent function $F(x)$ as
\begin{eqnarray*}
F(x) & = & (\T_6((g^3+g^2+1)x)+1)\T_6(gx^5)+\T_6((g^3+g^2+1)x)\T_6(g^{22}x^5) \\
& = & \T_6((g^3+g^2+1)x)\T_6((g+g^{22})x^5)+\T_6(gx^5).
\end{eqnarray*}
In view of Corollary \ref{simple}, $F$ is described as $F(x) = \T_6(gx^5)$ if 
$\T_6((g^3+g^2+1)x) = 0$ and  $F(x) = \T_6(g^{22}x^5)$ otherwise. \\
Expanding the trace terms we get the unique representation of $F(x)$ as a
polynomial of degree at most $2^6-1$, as one expects, as a rather complicated expression
\begin{align*}
F(x)&=g^{51}x^{56} + g^{27}x^{52} + g^{12}x^{50} + g^{39}x^{49} + g^2x^{48} +
g^3x^{44} + x^{42} + g^{54}x^{41} \\
&+ g^{24}x^{40} + g^{27}x^{38} + g^{24}x^{37} + g^{15}x^{35} + g^{33}x^{34} +
g^4x^{33} + g^7x^{32}+ g^{57}x^{28}\\
& + g^{45}x^{26} + g^6x^{25} + gx^{24} + g^{33}x^{22}
+ x^{21} + g^{12}x^{20} + g^{45}x^{19} + g^{48}x^{17}\\
& + g^{35}x^{16} +g^{60}x^{14} + g^{54}x^{13} + g^{32}x^{12} + g^{48}x^{11} + g^6x^{10}
+ g^{49}x^8 + g^{30}x^7 \\
& + g^{16}x^6 + g^3x^5 + g^{56}x^4 +g^8x^3 + g^{28}x^2 + g^{14}x.
\end{align*}
The Fourier spectrum of $F$ is the multiset $\{-8^{28},8^{36}\}$, where the integer 
in the exponent denotes the multiplicity of the corresponding Fourier coefficient 
in $\{\widehat{F}(b)\;|\;b\in \F_{p^n}\}$. The algebraic degree of $F$ is $3$.
\end{example}

We are particularly interested in a first direct construction of not weakly regular bent functions in
the framework of finite fields. Taking advantage of the binomial description of $2$-plateaued partially bent functions in 
Proposition \ref{binomial}, in the subsequent corollary we present an infinite class of not weakly regular bent functions in 
arbitrary odd characteristic with a simple description.
\begin{corollary}
For an odd integer $n$ divisible by $3$, let $\beta_1,\beta_2$ be solutions in $\F_{p^n}$ of $x^{p^2}+x^p+x$ (which are linearly
independent over $\F_p$), let $\delta \in \F_{p^n}$ such that $\Tr(\delta\beta_1^2) = l \ne 0$ and
$\Tr(\delta\beta_1\beta_2) = 0$ and let $\Gamma\in\F_p^n$ such that $\Tr(\delta\Gamma\beta_2) = t \ne 0$. 
Let $r=2^\kappa$, $\kappa\ge 1$, $a_0 = 1$, $c_0 = 0$, and for $a_k\in \F_p^*$ let $c_k\in\F_p$, $1\le k\le p-1$, be given as 
\begin{equation} 
\label{c}
c_k = t^{-1}[(1-a_k)\Delta+k],\;\mbox{where}\;\Delta = \Tr(\frac{p+1}{2}\beta_2^2+\beta_2^{p^r+1}).
\end{equation}
Then the function
\begin{equation}
\label{F,a} F(x) = a_{\Tr(\delta\Gamma x)}\Tr(\frac{p+1}{2}x^2+x^{p^r+1}) + c_{\Tr(\delta\Gamma x)}\Tr(\delta\Gamma x) 
\end{equation}
is bent. It is not weakly regular if and only if $a_k$ is a nonsquare in $\F_p$ for some $1\le k\le p-1$.
\end{corollary}
\bproof
By the proof of Proposition \ref{binomial}, the $2$-dimensional linear space of the partially bent function $f(x) = \Tr(\frac{p+1}{2}x^2+x^{p^r+1})$
consists of the solutions of $x^{p^2}+x^p+x$. Hence with this choice of $\beta_1,\beta_2$ and $\delta$ 
the requirements in Corollary \ref{simple} are satisfied. Let $f_k(x) = a_kf(x)$, then with the definition of
$c_k$ in $(\ref{c})$, the partially bent functions $f_k$ satisfy $f_k(\beta_2) + \Tr(\delta\Gamma_k\beta_2) = f_0(\beta_2) + k$,
i.e. the functions $g_k(x) := f_k(x) + c_k\Tr(\delta\Gamma x)$ satisfy $supp(\widehat{g_k})\cap supp(\widehat{g_j}) = \emptyset$
for $0 \le l\ne k\le p-1$. With Corollary \ref{simple} the function $F(x)$ is bent. \\
By Theorem 1 in \cite{aw}, since $n$ is odd, the nonzero Fourier coefficients of $f(x)$ change the sign if we multiply $f(x)$
by a nonsquare in $\F_p$ (see also \cite[Theorem 4.3]{agw}). Since we choose $g_0(x) = f_0(x) = f(x)$, we combine partially bent 
functions with Fourier coefficients of opposite sign if at least for one $1\le k\le p-1$ the coefficient $a_k$ is a nonsquare. 
For details on this technique of obtaining not weakly regular bent functions we refer to \cite{agw,aw}.
\eproof

We remark that $\beta_1,\beta_2$, the solutions of the linearized polynomial $x^{p^2}+x^p+x$ (which can be determined with standard methods
using linear systems) are independent of $n$.  For $x \in \F_{p^n}$, the value of the linear term  $k = \Tr(\delta\Gamma x)$ decides on the 
coefficient $a_k$ for the binomial part  in $(\ref{F,a})$  and the coefficient $c_k$ for the linear part in $F(x)$ in  $(\ref{c})$.


\section{Conclusion}
Based on earlier constructions of Boolean bent functions from near-bent functions
(see \cite{cpt,lg}), in \cite{agw,aw1} constructions of bent functions in arbitrary
characteristic have been presented.
Amongst others, with these constructions the first infinite classes of
not weakly regular bent functions were obtained. The functions
are given in multivariate form or as functions from
$\F_p^{n-1}\times \F_p$ to $\F_p$. Until then only sporadic examples
of not weakly regular bent functions have been found via computer search, 
\cite{hk,hk1,tyz}. All of these sporadic examples were given in univariate 
form, i.e. as functions from $\F_{p^n}$ to $\F_p$ (represented as polynomials 
in trace form). In this article an equivalent but more involved procedure in the 
framework of functions from $\F_{p^n}$ to $\F_p$ is developed.
In particular we obtain the first direct construction of not weakly regular 
bent functions in univariate form. We take advantage of some infinite classes of 
partially bent monomials and binomials with a $2$-dimensional linear space, to 
construct (not) weakly regular bent functions from $\F_{p^n}$ to $\F_p$ with a 
simple representation.

\section{Appendix}

\begin{example}
Let $g$ be a root of the primitive polynomial
$x^4+2x^3+2\in\F_3[x]$. Then by Lemma \ref{monomial}, the functions
\[ f_0(x)=f_2(x)=\T_4(g^4x^{28})\;\mbox{and}\;f_1(x)=\T_4(2g^4x^{28}) \]
from $\F_{3^4}$ to $\F_3$ are $2$-plateaued. For both functions the linear space is
$\Lambda=\la \beta_1,\beta_2\ra$ with $\beta_1=g^2, \beta_2=g^3+2g+1$.
We then have $\T_4(\beta_1^2)=2$, $\T_4(\beta_1\beta_2)=0$, $\T_4(\beta_2^2)\ne 0$.
With Proposition \ref{prop} we get the near-bent functions from $V_3$ to $\F_3$
\[ g_{0_{|3}}(x)=\T_4(g^4x^{28}),\; g_{1_{|3}}(x)=\T_4(2g^4x^{28}+x),\; g_{2_{|3}}(x)=\T_4(g^4x^{28}+2x) \]
with $supp(\widehat{g_j}_{|3})\cap supp(\widehat{g_k}_{|3})=\emptyset$ if $j\ne k$.
With Theorem \ref{theo}, we obtain the weakly regular bent function
\begin{align*}
F(x)&=2[(\T_4(\gamma x)-1)(\T_4(\gamma x)-2)g_0(x)+\T_4(\gamma x)(\T_4(\gamma
x)-2)g_1(x)\\
&+\T_4(\gamma x)(\T_4(\gamma x)-1)g_2(x)]
\end{align*}
or alternatively
\[
F(x) = \left\{\begin{array}{ll}
\T_4(g^4x^{28}) & \mbox{if}\quad\T_4(\gamma x) = 0 \\[.3em]
\T_4(2g^4x^{28}+x) & \mbox{if}\quad\T_4(\gamma x) = 1 \\[.3em]
\T_4(g^4x^{28}+2x) & \mbox{if}\quad\T_4(\gamma x) = 2,
\end{array}\right.
\]
with Fourier spectrum $\{(-9)^{21},
    (-9\epsilon_3)^{30},
    (-9\epsilon_3^2)^{30}\}$.
\end{example}

\begin{example}
\label{ex3} Let $g$ be a root of the primitive polynomial
$x^3+2x+1\in\F_3[x]$. According to Proposition \ref{binomial}, the
functions
\[ f_0(x)=\T_3(2x^2+x^{10}),\; f_1(x)=\T_3(2x^2+x^4) \]
are partially bent functions from $\F_{3^3}$ to $\F_3$. For each function, the linear space is
$\Lambda=\la \beta_1, \beta_2\ra$ with $\beta_1=g$ and $\beta_2=1$.
Again we have $\T_3(\beta_1^2)=2\ne 0$ and $\T_6(\beta_1\beta_2)=0$, but now $\T_6(\beta_2^2)=0$,
see Remark \ref{inepro}. \\[.3em]
{\bf A:} Put $f_2(x) = f_1(x)$, then with Proposition \ref{prop} we obtain the near-bent functions from $V_2$ to $\F_3$
\[ g_{0_{|2}}(x)=\T_3(2x^2+x^{10}),\; g_{1_{|2}}(x)=\T_3(2x^2+x^4+2g^2x),\; g_{2_{|2}}(x)=\T_3(2x^2+x^4+g^2x) \]
with $supp(\widehat{g_j}_{|2})\cap supp(\widehat{g_k}_{|2})
=\emptyset$  if $j\ne k$.
With $\gamma=2g$, Theorem \ref{theo} yields the weakly regular bent function
\begin{align*}
F(x)&=2[(\T_3(2gx)-1)(\T_3(2gx)-2)g_0(x)+\T_3(2gx)(\T_3(2gx)-2)g_1(x)\\
&+\T_3(2gx)(\T_3(2gx)-1)g_2(x)] \\
&= \left\{\begin{array}{ll}
\T_3(2x^2+x^{10}) & \mbox{if}\quad\T_3(gx) = 0 \\[.3em]
\T_3(2x^2+x^4+2g^2x) & \mbox{if}\quad\T_3(gx) = 2 \\[.3em]
\T_3(2x^2+x^4+g^2x) & \mbox{if}\quad\T_3(gx) = 1,
\end{array}\right.
\end{align*}
with the Fourier spectrum
$\{(-3^{3/2}i)^{9},(-3^{3/2}i\epsilon_3)^{12},(-3^{3/2}i\epsilon_3^2)^{6} \}$. \\[.3em]
{\bf B:} Put $f_2(x) = 2f_1(x) = \T_3(x^2+2x^4)$, such that the nonzero Fourier coefficients of $f_1$ and $f_2$ have
opposite signs. Again with Proposition \ref{prop}, we obtain the near-bent functions from $V_2$ to $\F_3$
\[ g_{0_{|2}}(x)=\T_3(2x^2+x^{10}),\; g_{1_{|2}}(x)=\T_3(2x^2+x^4+2g^2x),\; g_{2_{|2}}(x)=\T_3(x^2+2x^4+g^2x) \]
with $supp(\widehat{g_j}_{|2})\cap supp(\widehat{g_k}_{|2})=\emptyset$ if $j\ne k$, and then with Theorem \ref{theo} 
the not weakly regular bent function
\begin{align*}
F(x)&=2[(\T_3(2gx)-1)(\T_3(2gx)-2)g_0(x)+\T_3(2gx)(\T_3(2gx)-2)g_1(x)\\
&+\T_3(2gx)(\T_3(2gx)-1)g_2(x)] \\
&= \left\{\begin{array}{ll}
\T_3(2x^2+x^{10}) & \mbox{if}\quad\T_3(gx) = 0 \\[.3em]
\T_3(2x^2+x^4+2g^2x) & \mbox{if}\quad\T_3(gx) = 2 \\[.3em]
\T_3(x^2+2x^4+g^2x) & \mbox{if}\quad\T_3(gx) = 1,
\end{array}\right.
\end{align*}
where again $\gamma = 2g$. As one would expect, the unique polynomial representation of this not weakly regular bent
function does not look very simple:
\begin{align*}
F(x) &=g^8x^{24} + g^3x^{22} + x^{21} + x^{20} + gx^{19} + g^{17}x^{18} + g^{11}x^{16} +
    2x^{15} + g^{16}x^{14}\\
& + g^3x^{13} + g^8x^{12} + g^9x^{11} + g^5x^{10} +
    g^3x^9 + g^{20}x^8 + x^7 + g^{23}x^6 + 2x^5 \\
&+ g^{21}x^4 + gx^3 +
    g^9x^2 + 2x.
\end{align*}
The Fourier spectrum of $F$ is \\
$\{(3^{3/2}i)^{3}, (-3^{3/2}i)^{6}, (3^{3/2}i\epsilon_3)^{3},
(-3^{3/2}i\epsilon_3)^{9}, (3^{3/2}i\epsilon_3^2)^{3},
(-3^{3/2}i\epsilon_3^2)^{3} \}$. The algebraic degree of $F$ is $4$.
We remark that $F$ attains the upper bound on the algebraic degree
of bent functions, see \cite{aw1,h}.
\end{example}

\end{document}